\documentclass[11pt]{article}
\usepackage{graphicx}
\usepackage{amsmath}
\usepackage{amsfonts}
\usepackage{url}

\def\0{{\bf 0}}

\def\real{\mathbb{R}}

\newcommand{\enhe}{\~{n}}

\newcommand{\beq}{\begin{equation}}
\newcommand{\enq}{\end{equation}}
\newcommand{\beqa}{\begin{eqnarray}}
\newcommand{\enqa}{\end{eqnarray}}

\newcommand{\uno}{1\hspace{-2.8pt}\mbox{l}}
\newcommand{\neno}{\newline\noindent}

\newcommand{\esp}{\mathbb{E}}

\newcommand{\nn}{N\!N}

\newcommand{\X}{{\cal X}}
\newcommand{\Y}{{\cal Y}}
\newcommand{\cL}{{\cal L}}

\newcommand{\veno}{{\vskip 1pc \noindent}}

  \def\qed{\hfill\hbox{${\vcenter{\vbox{
        \hrule height 0.4pt\hbox{\vrule width 0.4pt height 6pt
        \kern5pt\vrule width 0.4pt}\hrule height 0.4pt}}}$}}
        
\usepackage{float}
\restylefloat{table}

\begin{document}

\begin{center}
\begin{Large}
Permutation tests in the two-sample problem for functional data
\end{Large}
\vskip 2pc
\begin{large}
A. Caba{\enhe}a\footnote{Departament de Matem\`atiques,
Facultat de Ci\`encies, 
Universitat Aut\`onoma de Barcelona. \\ Email: acabana@mat.uab.cat},
A. M. Estrada\footnote{Dpto. de Matem\'aticas,
Universidad de Los Andes, Bogot\'a. Email: am.estrada213@uniandes.edu.co}, 
J. I. Pe{\enhe}a\footnote{Dpto. de Matem\'aticas,
Universidad de Los Andes, Bogot\'a. Email: ji.pena45@uniandes.edu.co} and
A. J. Quiroz\footnote{Corresponding Author. Dpto. de
Matem\'aticas, Universidad de Los Andes. Address: Dpto. de
Matem\'aticas,
Universidad de Los Andes, Carrera 1, Nro. 18A-10, edificio H,
Bogot\'a, Colombia.
Phone: (571)3394949, ext. 2710. Fax: (571)3324427. Email: aj.quiroz1079@uniandes.edu.co} 
\end{large}
\end{center}
\vskip 1pc
\begin{small}
{\bf Abstract} Three different permutation test schemes are discussed and compared 
in the context of the two-sample problem for
functional data. One of the procedures was essentially introduced by Lopez-Pintado and Romo,
in \cite{lpr09}, using notions of functional data depth to adapt the ideas 
originally proposed by Liu and Singh,
in \cite{ls93} for multivariate data. 
Of the new methods introduced here, one is also based on functional data
depths, but uses a different way (inspired by Meta-Analysis)
to assess the significance of the depth differences. The second
new method presented here adapts, to the functional data setting, the $k$-nearest-neighbors
statistic of Schilling, \cite{sch86b}. The three methods are compared among them
and against the test of Horv\'ath and Kokoszka \cite{hk12} in simulated
examples and real data. The comparison considers the performance of the statistics in terms
of statistical power and in terms of computational cost.

\end{small}

\section{Introduction}
\subsection{The Two-Sample Problem for Functional Data}

Nowadays it is not uncommon to have access to data varying over time or space, sampled at a high 
enough rate, in such a way that it is valid to think that the whole ``curve'' or ``function'' of 
interest is available. Moreover, it is frequent the situation in which the data analyst has several
of these curves from different individuals (or repetitions from a fixed individual). In such 
circumstances, the data set can be thought as made of independent realizations of a function 
produced by some random mechanism. A typical example of this situation would involve the analysis of a set of heart beat signals (or some other relevant physiological electric signal) from different 
individuals. Another example would be a set of optical coherence tomography images of the retina and 
optic nerve in different individuals. Real data examples abound. The statistical tools that have 
been developed for this type of data form a body of methods called Functional Data Analysis (FDA). 
In this context, the statistician wants to answer the usual questions. For instance, what is a good  
estimate of center for the data set of curves? Should certain curves  be consider outliers, away 
from the main body of data? What are the principal ``directions'' of variation of the data set?, 
etc. 

Good reviews of the fundamental ideas available for Functional Data Analysis can be found in the 
books by Ramsay and Silverman \cite{rs1, rs2}, Horvath and Kokoszka \cite{hk12} and Ferraty and Vieu 
\cite{fv06}, where various  
application examples can be found. The subject of FDA can be thought of as a 
generalization of Multivariate Analysis to infinite dimension. Sometimes, the methods from 
Multivariate Analysis extend almost intact to the FDA context. Such is the case for Principal 
Component Analysis. In other instances, the methods proposed for FDA are far apart from their 
Multivariate Analysis counterparts and try to reflect the functional nature of the data. An 
important instance of this is depth measures. The depth measures that have been proposed for FDA in 
\cite{fm01}, \cite{lpr07} and \cite{lpr09} are clearly distinct from those reviewed in \cite{lps99} 
(see also \cite{ls93}) for finite dimensional multivariate data.

In the functional data setting, the two-sample problem can be stated as follows:
Let $X_1(t), \cdots$, $X_m(t)$ denote an i.i.d. sample of real valued curves defined
on some interval $J$. Let
$\mathcal L(X)$ be the common probability law of these curves. Likewise, let
$Y_1(t), \cdots , Y_n(t)$, be another i.i.d. sample of curves, independent
of the $X$ sample and also defined on $J$, with
probability law $\mathcal L(Y)$. We want to test
the null hypothesis, $H_0$: $\mathcal L(X) = \mathcal L(Y)$, against the general
alternative $\mathcal L(X) \neq \mathcal L(Y)$.

Fairly different approaches have been put forward in the literature for the
functional data two-sample problem. We describe next three of the relevant
ideas that have been proposed.

Hall and Van Kielegom \cite{hvk} consider
a bootstrap version of a generalized two-sample Cram\'er-von Mises statistic, based on
comparison of estimated probabilities of the events
\[
\mathcal L(X)\left(\{X\leq z\}\right)\qquad\mbox{and}\qquad 
\mathcal L(Y)\left(\{Y\leq z\}\right),
\]
where inequalities between functions are interpreted as holding when they
hold for every point. In this reference, the effect of smoothing on the
power of two-sample statistics is also discussed. The advice given by these
authors points in two directions: (i) Regarding smoothing, ``less'' should
generally be preferred over ``more'' and (ii) Whenever possible, the smoothing 
applied to the $X$ and $Y$ samples should be the same.
 
Horvath and Kokoszka (\cite{hk09} and \cite{hk12}) consider a statistic
which is a quadratic form based on the (Hilbert space) inner products of
the curves from both samples with the  principal components estimated from the covariance operator of
the pooled sample. Two versions of this test are considered in \cite{hk12}, one of which we will
include in our comparisons below, since both from the practical implementation
and performance viewpoints, it seems to us a very relevant procedure for the two-sample problem.  
In connection with this type of statistic, see also Benko et al \cite{benko}.

Mu{\~n}oz Maldonado et al \cite{munoz}, 
in the context of studying the physical process
of aging of the brain, consider the two-sample problem on registered curves
of tissue density profiles of the brain in old versus young rats.
As similarity measure between curves, they use Pearson's correlation coefficient
applied to the vectors obtained by registering the curves on a common grid. 
From the within groups correlation values, three different statistics are
considered. Relevant $p$-values for these statistics are obtained
by means of the natural permutation procedure, that we outline
for the reader's convenience: As above, let $Z$ denote the
pooled, or joint, sample, made of the union of the $X$ and $Y$ samples.
\veno{\it Procedure $p$-value}\label{ppval}

\noindent (i){\it
From the joint sample, $Z$, select a random subset
of size $m$. Declare that the chosen elements belong to the $X$ sample and the remaining $n$ to the
$Y$ sample}. 
\neno (ii) {\it Compute the statistic of interest for this artificial pair of samples.}
\neno (iii){\it Repeat steps (i) and (ii) a large number $B$ of times, ($B$ = 10,000,
for instance) and}
\neno (iv) {\it From the $B$ values computed above, extract an approximate $p$-value
for the observed statistic (the one calculated with the original $X$ and $Y$ samples)}  

The same basic permutation procedure for $p$-value estimation that we have just outlined, 
is used in the methods to be presented and
evaluated in the present article. These methods differ significantly in nature from the ones considered in \cite{munoz}.

\subsection{The Permutation Test Principle and Schilling's Statistic}

The {\it Permutation Test Principle} can be found behind a good portion of Non-Parametric 
Statistical methodology.
For instance, the Mann-Whitney-Wilcoxon statistic for the two-sample problem (in the case of 
difference in location) and the Ansari-Bradley statistic for the difference in scale two-sample 
problem can be formulated as permutation tests 
and their theory derived through this viewpoint. The same holds for the Kruskal-Wallis One-Way-Anova 
and for the Non-Parametric correlation measures of Spearman and Kendall. Complete discussions of the 
theory of all these procedures can be found in the classical text of Randles and Wolfe \cite{rw79}. 
Permutation tests are also behind more modern 
procedures, as the graph-theoretic generalization to the multivariate setting of the Wald-Wolfowitz 
two-sample test 
proposed by Friedman and Rafsky \cite{fr79} and the $k$-nearest-neighbor multivariate two-sample test of 
Schilling \cite{sch86b}. Since this last one will be adapted here to the Functional Data context, we 
will describe it below in some detail.

Regarding the use of permutation procedures, in his book \cite{good05}, Phillip Good claims:
\begin{quote}
``Distribution-free permutation procedures are
the primary method for testing hypotheses. Parametric procedures and the
bootstrap are to be reserved for the few situations in which they may be
applicable.''
\end{quote}

Good explains some of the reasons for the convenience of using permutation tests:
\begin{quote}\begin{itemize}
\item[]``Pressure from regulatory agencies for the use of methods that yield exact
significance levels, not approximations.''
\item[] ``A growing recognition that most real-world data are drawn from mixtures
of populations.''

\end{itemize}
\end{quote}

On one hand, the views quoted above reflect the convenience of considering non-paramet\-ric procedures
in the context of ``Big-Data'', in which the validity of a model for the huge data set at hand is fairly unlikely
to hold (the most one could hope for in that situation is to have a population coming from a mixture of 
distributions). On the other hand, the push towards the employment of permutation procedures reflects the increasing 
confidence of the statistical community on the availability of powerful computational resources.
In addition, the theoretical power analysis carried out by Randles and Wolfe \cite{rw79} in important examples,
shows that in many relevant cases, permutation methods do not lose significantly,
in terms of power, with respect to their optimal parametric counterparts, when the parametric
assumptions do hold.

We now describe, borrowing in part from \cite{aq06}, the test of Schilling \cite{sch86b} in some detail. Next, 
we will explain how, when this test is viewed as a permutation test, its 
null distribution can be approximated very efficiently, by a ``table permutation'' algorithm. The main 
features of Schilling test, as described in the present section, will essentially remain unaltered 
in the Functional Data version.

Suppose we have samples $X_1,\dots,X_m$ i.i.d. from a distribution $P_X$ and  $Y_1,\dots,Y_n$ i.i.d. from
a distribution $P_Y$, both living 
in $\real^d$. To test the null hypothesis $P_X=P_Y$ against the alternative $P_X\neq P_Y$, the 
method proposed in \cite{sch86b} is as follows: Let $N=m+n$. Denote by $Z_1,\dots,Z_{N}$ 
the pooled sample obtained by concatenation of the $X$ and $Y$ samples. Fix an integer $k>0$.
\veno{\it Schilling's Procedure}\label{schproc}
\begin{itemize}
\item[1.] For each $i\leq N$, find the $k$-nearest-neighbors, with respect to Euclidean distance, of 
$Z_i$ in the pooled sample. Let $N\!N_i(r)$ represent the $r$th nearest neighbor to $Z_i$, for $r
\leq k$. Assuming (under the null hypothesis) that the common underlying probability distribution is 
continuous, all the nearest neighbors are uniquely defined with probability 1. Otherwise, break 
distance ties at random. 
\item[2.] For each $i\leq N$ and $r\leq k$, compute the indicator variables   
\[   
\begin{array}{cccl}   
I_i(r)&= &1 &\mbox{ if } N\!N_i(r) \mbox{ belongs to the same sample as } Z_i \\   
&= &0 & \mbox{ otherwise. }   
\end{array}   
\]
\item[3.]  
Compute the statistic    
\begin{equation}\label{ecu1}  
T_{N,k}=\frac{1}{Nk}\sum_{i=1}^N\sum_{r=1}^kI_i(r).   
\end{equation}  
\end{itemize}
$T_{N,k}$ is the proportion of all neighbor pairs in which a point and its neighbor belong to the 
same sample. The rationale for considering $T_{N,k}$ is simple: when $H_0$ does not hold, points from the same sample will tend to clump together in those regions where their probability density is larger than the other, causing a high value of the statistic.   

In a separate paper, Schilling \cite{sch86b} provides the limiting theory for his statistic 
$T_{N,k}$. The main assumptions are that the common distribution has a density   
$f$ that is continuous on its support, and that $m/N$ and $n/N$ converge to non-zero   
limits $\lambda_1$ and $\lambda_2$, respectively.   
The limiting distribution depends on the probability, under the null assumption that
the samples come from the same distribution, that $Z_1$ and $Z_2$ (or any two 
vectors of the pooled sampled) are mutual   
nearest neighbors, as well as the probability of $Z_1$ and $Z_2$ sharing   
a nearest neighbor. More precisely, if for each $r$ and $s$ in $\{1,\dots,k\}$, one lets   
\begin{itemize}   
\item[(i)] $\gamma_N(r,s)=\Pr(N\!N_1(r)=Z_2,\>N\!N_2(s)=Z_1)$ \hspace{5pt} and   
\item[(ii)] $\beta_N(r,s)=\Pr\left(\nn_1(r)=\nn_2(s)\right)$, for $1\leq r,s\leq k$,   
\end{itemize}   
it turns out that $N\gamma_N(r,s)$ and $N\beta_N(r,s)$ have positive limits, as $N\rightarrow   
\infty$, which do not depend on the underlying $f$. These limits, together with $\lambda_1$ and 
$\lambda_2$, determine the asymptotic distribution of $T_{N,k}$. See \cite{sch86b} for more details.   

Let us now turn the discussion towards the permutation test nature of $T_{N,k}$ and its efficient implementation as such.

Conditioning on the observed points of the pooled sample, under the null hypothesis, the labelling of any of these
points as $X$ or $Y$ is, essentially, an arbitrary choice of the experimenter, since all elements in the
pooled sample were generated by the same random mechanism. Thus, the permutation distribution
is a valid reference for the observed (true) value of the statistic considered. It is easy
to see that for this permutation distribution, the expected value of $T_{N,k}$ is
\[
\esp T_{N,k}=\esp I_i(r)=\frac{m(m-1)+n(n-1)}{N(N-1)},
\]
while the variance depends on the amounts of pairs of points that are mutual neighbors and the amount of
pairs of points that share a common neighbor. As explained above, under the alternative, we expect $T_{N,k}$ 
to take values above the null expected value. 

Now, conditionally on the observed pooled sample, the set of $k$ nearest neighbors of each of the pooled sample points
is completely determined. In order to sample from the null permutation distribution of $T_{N,k}$, the $k$ nearest 
neighbors of every point in the pooled sample must be identified. This is the heaviest computational burden of the
procedure, but it needs to be performed only once!, as we now explain. Suppose that, in the pooled sample, we 
originally keep the natural ordering from concatenation of the samples: 
$Z_1$ is $X_1$, $Z_2$ is $X_2$,$\dots$, $Z_m$ is $X_m$, $Z_{m+1}$ is $Y_1$, $Z_{m+2}$ is $Y_2$, $\dots$, 
$Z_N$ is $Y_n$. When the nearest neighbors are computed, a $N\times k$ table is constructed that in the 
$i$-th row contains the indices of the $k$ elements nearest to $Z_i$ in the pooled sample. In order to compute
$T_{N,k}$ it is enough to count, in each row in the first $m$, how many of the indices are less than or equal to $m$, and in each row of the last $n$, how many of the indices are greater than $m$. 
For an iteration of the permutation procedure, a random permutation $\sigma$ is applied to the indices 1,2,
$\dots$,$N$, and, to compute $T_{N,k}$, in each row of the original nearest neighbors table we must count the
number of indices that belong in $X^{\sigma}=\sigma\left(\{1,2,\dots,m\}\right)$ or in its complement. 
Now, an index $l$ belongs
to $X^{\sigma}$ if, and only if  $i=\tau(l)$ for some $i\leq m$ and
$\tau=\sigma^{-1}$. Since the distribution of $\sigma$ is Uniform in the set of permutations, so is the 
distribution of $\tau$. It follows that, for an iteration of the random permutation procedure, the following
suffices: 
\neno (i) To the elements of the original nearest neighbor table, apply a random permutation $\tau$.
\neno (ii) In the $r$-th row, if $\tau(r)\leq m$, count the number of elements that are less than or equal to $m$, 
else, count the number of elements greater than $m$. 
 \neno (iii) with the numbers obtained in (ii) compute $T_{N,k}$. \neno
 All these operations (including the generation of $\tau$) can be performed in time O$(Nk)$, thus the
 computational cost of an iteration is basically linear in $N$. As for the initial cost of setting
 the neighbors indices table, several sub-quadratic algorithms have been developed for this problem,
 since the fundamental contribution of Friedman et al \cite{fbs75}. Although those algorithms are intended
 to be used on data living in Euclidean space, they still work fine in relatively large dimensions
 (in the order of a few hundred coordinates), and therefore can be adapted to functional data when the
 curves have been registered on a common grid with no more than hundreds of points.   

\subsection{Depth Measures for Functional Data} 
In this subsection we describe two relevant measures of data depth for functional data. These measures
seek to provide a measure of the outwardness of a curve with respect to a functional data set. Functional data depth
can be used to rank the curves in a data set, from center outward.

In order to present the functional data depth of Fraiman and Muniz \cite{fm01} consider first a univariate sample, 
$U_1,\dots,U_n$, let $U_{(1)},\dots,U_{(n)}$
denote the corresponding order statistics. Assume these are uniquely
defined (no ties). Then, if $U_i=U_{(j)}$, the natural depth of $U_i$ in the sample 
is given by 
\beq\label{ecu2}
D_n(U_i)=\frac{1}{2}-\left|\frac{1}{2}-\left(\frac{j}{n}-\frac{1}{2n}\right)\right|.
\enq
This notion of depth assigns minimal and equal depth to the two extreme
values of the sample, maximum depth to the innermost point (or points in case
$n$ is even) and changes linearly with the position the datum occupies in the sample. 

In a sample of functional data, ${\cal X}=\{X_1,X_2,\dots,X_n\}$, defined
on a common interval $J$, for each fixed
$t$ in $J$ we compute the univariate depth defined above of each value $X_i(t)$,
with respect to the sample $X_1(t),X_2(t),$ $\dots,X_n(t)$. Call this 
depth $D_n(X_i(t))$. Then, the depth of the $i$-th curve with respect to the
$X$ sample, $\cal X$, is given by
\beq\label{ecu4}
I(X_i,{\cal X})=\int_JD_n(X_i(t))\mbox{d}t.
\enq
In practice, the integral of the definition is approximated after computing
the univariate depths on a finite grid of values of $t$.

Other functional data depths, considered by Lopez-Pintado and Romo in \cite{lpr07} and \cite{lpr09}, are based on the
notion of a functional band. Given real functions $u_1,\dots, u_r$, defined on an
interval $J$, the band defined by these functions is the two dimensional set
\[
V(u_1,\dots, u_r)=\{(t,y)\in\mathbb{R}^2:\> t\in J\mbox{ and } 
\min_{j \leq r}u_j(t)\leq y\leq \max_{j \leq r}u_j(t)\}.
\]
$V(u_1,\dots, u_r)$ is the region between the pointwise minimum and maximum of the
functions considered. For a function $u$ defined on the interval $J$, let its graph
be defined in the usual way: $G(u)=\{(t,u(t)): t\in J\}$. Then, the $r$-th band depth of a
function $u$ with respect to a functional data set $u_1,\dots, u_n$ is given by
\beq\label{ecu3}
S_n^{(r)}(u)=\binom{n}{r}^{-1}\sum_{1\leq i_1<\dots <i_r\leq n}
\uno_{G(u)\subset V(u_{i_1},\dots, u_{i_r})},
\enq
where, for an event $E$, $\uno_E$ takes the value 1 if $E$ occurs and 0 otherwise. This definition says
that the depth of a curve with respect to a sample is the fraction of all possible bands (of $r$ functions
from the sample) that contain the given curve. This definition tends to produce very low depth for most 
curves in the sample. In order to improve its performance, it was modified by L\'opez-Pintado and Romo
\cite{lpr09}, changing the indicator function $\uno_{G(u)\subset V(u_{i_1},\dots, u_{i_r})}$ by a
measure of the percentage of time that the curve $u$ remains in the band $V(u_{i_1},\dots, u_{i_r})$. 
See more details in \cite{lpr09}. 

The depth measure $I(X_i,{\cal X})$ of Fraiman and Muniz (\ref{ecu4}) has been shown to have good behaviour
in many contexts and is relatively easy to compute. For these reasons, we will use it in the
comparisons to be described in Section 3.  
 
\section{Three ways of implementing permutation tests for Functional Data}

\subsection{Wilcoxon's statistic}

In the context of multivariate data, Liu and Singh, \cite{ls93}, propose the consideration of the Wilcoxon
rank test, using multivariate depth measures, such as Tukey's half space depth or Liu's simplicial depth,
for instance, instead of the univariate ranks of the original Wilcoxon statistic. This is a natural proposal,
since one of the uses of multivariate depth measures is to provide a center-outward ranking of multivariate data. 
Briefly, if the depths of all points in the pooled sample are computed (with respect to that joint sample),
then these depths can be ordered, say from larger to smaller, and then a rank in $\{1,2,\dots,m+n\}$ can be
assigned to each point according to its depth. For definiteness, we can establish that rank 1 is given to
the data point of maximum depth and rank $m+n$ goes to the (outermost) point of minimal depth. Ties can be
resolved in the usual manners (see for instance \cite{hwc}). An important advantage of using Wilcoxon
method is that the reference quantiles depend on the sum of random subsets of $\{1,2,\dots,m+n\}$,
or a similar set, in case of ties. More recently, it has been suggested in \cite{lpr07} that a similar
adaptation of Wilcoxon's statistic can be made in the context of functional data, using an appropriate
notion of functional data depth. This statistic, using the depth of Fraiman and Muniz for functional
data, described in the
previous section, will be the first statistic included in our comparison.

\subsection{Combining depths through Meta Analysis}

Let again $\X =\{X_1, \cdots$, $X_m\}$ denote our functional $X$ sample and 
$\Y=\{Y_1, \cdots , Y_n\}$ the functional $Y$ sample. For each $X_i\in\X$, we consider
its depth with respect to the $Y$ sample with $X_i$ added. We denote this depth
$I(X_i,{\Y}\cup\{X_i\})$, following the notation in (\ref{ecu4}). 
This is a measure of how outlying the curve $X_i$ is with respect to the $Y$ sample.
If ``many'' of the $X_i$ turn out to be outlying with respect to $\Y$, that would
be evidence against the null hypothesis of equality of distributions. Similarly
we can measure how outlying is each curve $Y_j$ with respect to the $X$-sample, $\X$,
by computing $I(Y_j,{\X}\cup\{Y_j\})$. The first question is how to combine the values
of $I(X_i,{\Y}\cup\{X_i\})$, for all $i\leq m$, in a single number that combines
the information in all these depths. For this purpose, we rely in an idea coming
from Meta-Analysis.

To the depth $I(X_i,{\Y}\cup\{X_i\})$ we associate an empirical $p$-value,
\beq\label{ecu5}
p_i=\frac{\#\{j\leq n:I(Y_j,\Y\cup\{X_i\})\leq I(X_i,{\Y}\cup\{X_i\})\}}{n+1},
\enq
where $\#$ stands for cardinality (of a finite set). A small value of $p_i$
corresponds to an outlying $X_i$, in terms of depth. Since the $X_i$
are i.i.d., we can think of the values $\{p_1,\dots,p_m\}$ as a set of
nearly independent empirical $p$-values based on the depths of the $X_i$
with respect to the $Y$ sample. Actually, ignoring ties in the depth values,
the distribution of each $p_i$ under the null distribution is the discrete
uniform distribution on the set $\{1/n,2/n,\dots,1\}$.

In Meta-Analysis, the problem of combining the $p$-values for independent
tests of the same null hypothesis has been considered. One of the methods discussed
in the classic text of Hedges and Olkin, \cite{ho85}, is the following:
Reject the null hypothesis for large values of 
\beq\label{ecu6}
S_X=-\sum_{i=1}^m\ln p_i,
\enq
where the sum goes from 1 to $m$ in our context, since this is the number
of $p$-values to combine. The sub-index $X$ refers to the fact that we are
computing $p$-values for the $X$ curves. The rational 
for using $S_X$ is the following: In the continuous case, each $p_i$ would
have a Uniform distribution on $[0,1]$. Thus, each $-\ln p_i$ is an exp(1)
variable and, assuming independence of the $p_i$, 
$S_X$ will have a Gamma distribution, with
shape parameter $m$. Furthermore, large values of $S_X$ correspond to several
of the $p_i$ being small (close to zero), which is what we are interested
in detecting.

In our case, independence of the $p_i$ cannot be postulated, since all of
them are computed with respect to the same $Y$ sample. Preliminary evaluations
show that, for values of $m$ and $n$ in the few hundreds, the approximation
of the distribution of $S_X$ to the corresponding Gamma distribution,
is not satisfactory in our context. 
Still, we can use the statistic $S_X$ in a permutation procedure. Actually,
we describe next two manners of evaluating $p$-values based on $S_X$ and its
symmetric counterpart, $S_Y$.

In order to have a symmetric statistic, to the depths, $I(Y_j,{\X}\cup\{Y_j\})$,
of the $Y_j$ with respect to the $X$ sample, we associate the corresponding
empirical $p$-values
\beq\label{ecu7}
q_j=\frac{\#\{i\leq m:I(X_i,\X\cup\{Y_j\})\leq I(Y_j,{\X}\cup\{Y_j\})\}}{m+1},
\enq
and with these, compute the statistic $S_Y=-\sum_{j=1}^n\ln q_j$, for the depths
of the $Y_j$ respect to the $X$ sample. As statistic, we use the maximum of
$S_X$ and $S_Y$, $S=\max(S_X,S_Y)$. The reason for considering the maximum is the following:
When the two samples display a difference in ```scale''', it can happen that all, or most,
of the curves of the $X$ sample, turn out to be central with respect to the $Y$ sample
and $S_X$ will not show a significant value. In such a situation, typically, 
several curves of the $Y$ sample will turn out to be clearly outlying respect to the
$X$ sample, and the maximum will reach a significant value. 

In order to associate significance to the observed value of $S$, we apply the
permutation ``procedure $p$-value'' described at the end of subsection 1.1. Namely,
in each iteration, a random subset of size $m$ is chosen from the joint sample of
functional data and labelled as the $X$ sample, while the remaining set of curves
is labelled as the $Y$ sample. On these samples the statistics $S_{X}$, $S_Y$ and
$S$ are computed. From a large number of repetitions, the $p$-value of $S$ can be estimated.
The procedure just described is called MA1 (for Meta Analysis 1) 
in the example evaluation section (Section 3).

A second way of associating a $p$-value to the pair ($S_X$,$S_Y$) is based on the following:
\veno{\it Lemma: combining $p$-values}\label{lema}

{\it Let $p_X$ ($p_Y$) denote the $p$-value of $S_X$ ($S_Y$), under the null permutation 
distribution, as obtained from procedure $p$-value if all subsets of size $m$ were used (instead
of just a sample of size $B$) and assuming the null hypothesis. Then
\neno (i) $\Pr(p_X\leq t)\leq t$ for any $t\in (0,1)$, and the same holds for $p_Y$.
\neno (ii) $\Pr(2\min(p_X,p_Y)\leq t)\leq t$ for any $t\in (0,1)$. 
}  
\veno{\it Proof:} The null permutation distribution of $S_X$ 
is a discrete distribution that can not be assumed uniform on its range (some values
of $S_X$ can appear more frequently than others when subsets are chosen at random). This
is why part (i) of the Lemma is not obvious. Let $F$ denote the null permutation c.d.f.
of $S_X$ and let $S_{X,\mbox{\tiny obs}}$ denote the observed value of $S_X$. Recall that
large values of $S_X$ are considered significant. Then, 
clearly, $p_X=1-F(S_{X,\mbox{\tiny obs}}^-)$ and, for $t\in (0,1)$,
\[
\Pr(p_X\leq t)=\Pr(F(S_{X,\mbox{\tiny obs}}^-)\geq 1-t)=\sum_{\{s:F(s)> 1-t\}}\Pr(S_X=s)\leq t,
\]
by definition of $F$.

Since $p_X$ and $p_Y$ are not independent, to prove (ii) we can use (i) together with the usual union bound:
\[
\Pr(\min(p_X,p_Y)\leq t/2)\leq \Pr(p_X\leq t/2)+\Pr(p_Y\leq t/2)\leq t/2+t/2=t.
\]
\veno
Part (ii) of the Lemma tells us that an appropriate $p$-value for the ``statistic'' 
$2\min(p_X,p_Y)$ is the observed value of this statistic itself. Thus, our second way of getting
a $p$-value from $S_X$ and $S_Y$ is the following: Compute, approximately, $p_X$ and $p_Y$
for $S_X$ and $S_Y$, respectively, using the procedure $p$-value described above and  
use $2\min(p_X,p_Y)$ as $p$-value. In the experiments described in Section 3, $p_X$ and $p_Y$
will be computed using independent ``draws'' of subsets of the joint sample, a procedure
that yields good power results. This way of computing $p$-values is called MA2 in the
evaluations in Section 3.

\subsection{Schilling's Statistic}
The third possibility of permutation test for the two-sample problem for functional data,
considered here,
is Schilling's statistic, described in detail, for the multivariate setting, in the previous section.

In the case of functional data, the first step needed to set up Schilling's procedure is
the construction of the inter-curve distance matrix, $D=(d_{i,j})_{i,j\leq N}$, where $N=m+n$, is the size
of the joint sample. For this purpose, it seems natural to use the $\cL^2$ distance whenever possible.
In practice, if the functions have been registered on a common grid, say
$0=t_0<t_1<t_2<\dots <t_L=T$, a reasonable approximation to the distance between functions $Z_i$ and $Z_j$,
would be the $\cL^2$ distance based on the points of the grid:
\beq\label{ecu8}
d_{i,j}=\sum_{l=1}^L \Delta_l(Z_i(t_l)-Z_j(t_l))^2, \>\>\mbox{where}\>\> \Delta_l=t_l-t_{l-1}, \>\>
\mbox{for}\>\> l=1,2,\dots L.
\enq
We use this approximation in the Monte Carlo example considered in Section 3. If the grid used is equally
spaced, the $\Delta_l$ term in the sum in (\ref{ecu8}) can be omitted and 
the curves can be treated as points in $\real^L$ in order to compute faster the $k$-nearest-neighbors
of each data point by means of the algorithms developed for Euclidean data as described in
\cite{fbs75}, \cite{amnsw98} or \cite{bkl06}, for instance. When no common grid is available, the
distance matrix can be calculated after the functions in the joint sample
have been represented in terms of local polynomials, or some other basis functions, and
the $k$-nearest-neighbors identified by a quadratic algorithm (in the joint sample size $N$).
In the evaluations in Section 3, Schilling's procedure is implemented with number of neighbors $k=5$
and $k=10$.

The following section describes the comparison of these three permutation procedures among them and
against the principal components method of Horv\'ath and Kokoszka (\cite{hk09} and \cite{hk12})
that uses the scalar products of the functions in the two samples against common principal
components to produce a statistic with an approximate chi-squared distribution.

\section{Empirical comparison of the tests}

We compare the performance of our tests in a series of simulated experiments, and also test them in real data.

In order to enrich the comparisons, we have also computed the empirical power for Horvath and 
Kokoszka's principal component test for equality of mean functions, \cite{hk12}, which we consider
one of the best tests available for the two-sample problem in the functional data context. In the Horvath and 
Kokoszka's test, the null hypothesis that the mean functions of the functional samples $X_1, \dots,X_m$ and $Y_1, \dots,Y_n$ are equal is 
rejected for large values of the the statistic $U_{m,n}=\frac{nm}{n+m}\int_0^1(\bar X_m(t)-\bar Y_n(t))^2dt$, or more 
precisely, a projection version of $U_{m,n}$ that uses only the first $d$ terms in the $L^2$ expansion of $(\bar 
X_m(t)-\bar Y_n(t))$ in terms of the eigenvalues of the empirical covariance function
\begin{small}
\[
\hat z_{m,n}(t,s)=\frac{n}{(m+n)m}\sum_{i=1}^m(X_i(t)-\bar X_m(t))(X_i(s)-\bar X_m(s))+\frac{m}{(m+n)n}\sum_{i=1}^n(Y_i(t)-\bar Y_n(t))(Y_i(s)-\bar Y_n(s)).\nonumber
\]
\end{small}
\subsection{A simulation experiment}

We have simulated thirteen samples of functional data as realizations from a geometric Brownian motion process $f(t) = 
X_0 \exp{(rt - \frac{t\sigma^2}{2} + \sigma B_t)}$, where $r$ and $\sigma$ are, respectively,  the trend (drift) and 
volatility coefficients, $B_t$ is a standard Wiener process and  $X_0$ is the initial value.

Each function $f(t)$ is defined for $t \in [0, 2]$. More precisely, we made a partition in [0, 2] of 601 points, so $t 
= i/300$, $i = 0, ..., 600$, and we simulated the corresponding values for $f(t)$.
 
The first sample, which will be called $X$ from now on, consists of $m = 250$ realizations of a geometric Brownian 
motion with $\sigma = r = X_0 = 1$. The next four samples, denoted by $Y_{X_{1.25}}$, $Y_{X_{1.5}}$, $Y_{X_{1.75}}$ 
and $Y_{X_2}$, consist of $n = 200$ realizations of the same process with $\sigma = r = 1$ and $X_0 = 1.25$, $1.5$, 
$1.75$, $2$, respectively.
The last eight samples, $Y_{r_{1.25}}$, $Y_{r_{1.5}}$, $Y_{r_{1.75}}$, $Y_{r_2}$, $Y_{\sigma_{1.25}}$, 
$Y_{\sigma_{1.5}}$, $Y_{\sigma_{1.75}}$ and $Y_{\sigma_2}$ are defined accordingly, changing only one of the 
parameters at a time, and leaving the others constant at $1$. We could think of this Y samples as contaminated samples 
because of their different level, trend or volatility.

We want to test the null hypothesis $H_0$: $\mathcal L(X) = \mathcal L(Y)$, against the  alternative $\mathcal L(X) 
\neq \mathcal L(Y)$, where $Y$ is any of the `contaminated' samples previously defined. We also want to test the null 
hypothesis for two different reference samples ($\sigma = r = X_0 = 1$) with different sample sizes ($m = 250$, $n = 
200$) to assess the type I error of the procedures.

Some other relevant details of our implementation of the tests are as follows:
Schilling's method was
implemented with $k=5$ and $k=10$ neighbors. Larger values of $k$ were also considered, but did not
produce a clear improvement in power. For the Wilcoxon rank based procedure, we used a random tie breaking scheme and 
the standard normal 
approximation since the sizes of the samples were large enough. For the permutation tests (other than
Wilcoxon's) we used the permutation procedure described in section 1.2 with $1000$ replications to approximate the (conditional) permutation distribution. 
All permutation tests were replicated $N = 
100$ times, in order to estimate the power. Horvath and Kokoszka's test was performed  $1000$ times using the projections over the first four principal components of the covariance operator. 

\subsection{Results}

The objective was to evaluate the statistical power and the computational cost of every method. In order to compare 
the statistical power we counted how many times in the 100 repetitions of the experiments did the test reject
the null hypothesis when the alternative hypothesis was true. As for the analysis of computational cost, we used a 
desk top computer with an Intel Core i7 processor of 2,00 GHz and a RAM memory of 8 GB. 

The results obtained are summarized in Table 1 for a theoretical level of 5\%. The results for levels 1\% and 10\% were also computed, but are not included,
to save space, since they do not show an essentially different behaviour.

In Table \ref{tests} it is evident that Wilcoxon's statistic performs very well against volatility variations, even for small changes in the volatility parameter. But this statistic fails noticeably for
the other alternatives considered in our experiment. On the other hand, Horvath and Kokoszka's test (HK), being a test conceived for changes in the mean, shows the best performance against changes in the drift parameter, while its power numbers against changes in the origin (initial level) are good too, in general, although not among the top three. But HK results ineffective in picking the volatility changes.

The Meta Analysis methods have a power similar to HK against changes in the origin, while their power, although
reasonable, is inferior to HK's when it comes to changes in drift. On the other hand, both Meta Analysis procedures 
display excellent power against the volatility alternatives, where HK fails. Columns MA1 and MA2 in Table 1, are the
first columns that show appreciable power against all alternatives. In all cases, MA2 performs better than MA1 for 
small deviations from the null.    

Schilling's statistic (with $k=5$ and $k=10$), shows very good power against all the alternatives considered in our 
experiment. Overall, Schilling's statistic displays the best performance in terms of power among the methods
evaluated. The power figures for $k=10$ are the best of all tests included in our analysis.  

As for computing time, Wilcoxon's  and Horvath and Kokoszka's methods are the fastest, taking about
a second to produce the $p$-value in one of our geometric Brownian motion experiments. 
Meta Analysis methods are the slowest, taking around six minutes in order to obtain an approximate
p-value. In this respect, Schilling's procedure has an intermediate behaviour, taking some 
25 seconds for each approximate $p$-value to be obtained. 

\subsection{A real data example}

We compared the three permutation methods considered on real data sets drawn from hourly measurements of nitrogen 
dioxide (NO$_2$) in four neighbourhoods in Barcelona, Spain, namely Sants, Palau Reial, Eixample and Poblenou.  
Nitrogen dioxide, a known pollutant, is formed in most combustion processes using air as the oxidant. The measurements 
were taken along the years 2014 and 2015 in automatic monitoring stations and were obtained from \url{http:// 
dtes.gencat.cat/icqa}.  Contamination levels, a priori, could be different during working and non-working days, hence 
we partitioned the data accordingly: one set made out of hourly measures for all working days (220 functional 
observations approximately, each with 24 registrations) and another set made out of hourly measures in non-working 
days (120 functions approximately) each year. Other questions of interest are whether the levels of NO$_2$ changed 
from one year to the next in each neighbourhood, and the comparison of the pollution levels among the different 
neighbourhoods. In the following analysis, HK test was implemented using five principal components.

\begin{figure}[h]
\caption{Levels of NO$_2$ in non-working days in Sants and Palau Reial, 2014}\label{santspalau}
\begin{center}
\includegraphics[scale=0.5]{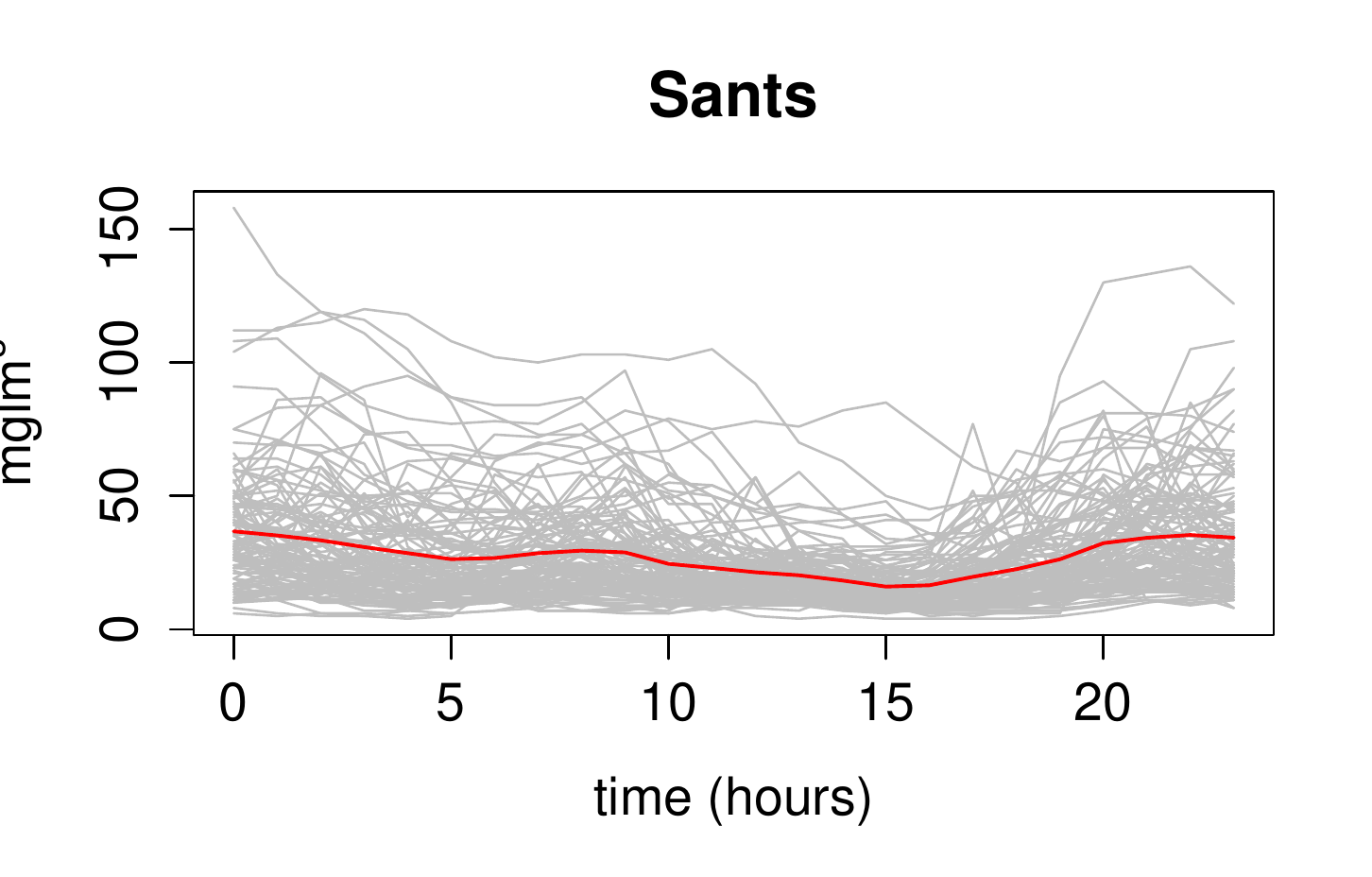}
\includegraphics[scale=0.5]{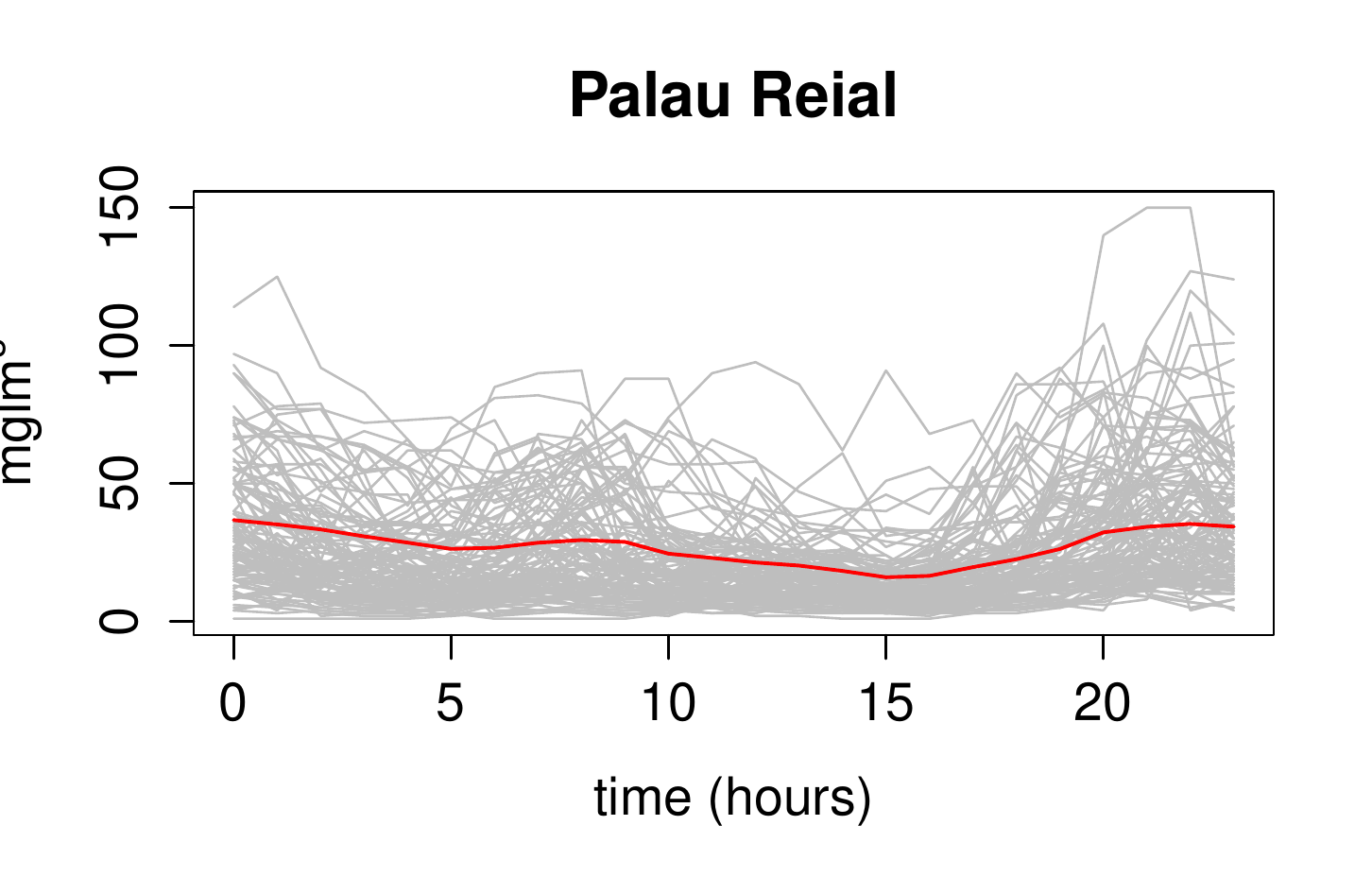}
\end{center}
\end{figure}

Figure \ref{santspalau} shows the levels of pollutants in non-working days in Sants and Palau Reial in 2014 in gray, ant the respective pointwise mean functions in red.

The next table shows the $p$-values obtained in some of the possible comparisons. We include only one of
the Meta Analysis procedures (MA2) and Schilling's statistic only for $k=10$.
The first three comparisons in Table \ref{compcont} correspond to comparing different neighbourhoods keeping
the year and type of day fixed. In this case, the permutation tests based on depth, Wilcoxon and MA2, are the
ones that work best, finding strong evidence of difference in all cases, while Schilling's test and the 
HK statistic fail to detect the differences or find only marginal evidence in some cases. In the next two comparisons, a fixed neighbourhood is compared against itself in the following year, keeping the type of day fixed. In these cases, all the methods reach the same conclusion: From 2014 to 2015, the pollution in Sants did not change
noticeably on non-working days, but significant changes are found from one year
to the next, in this neighbourhood, on working days, with
the Wilcoxon and Schilling procedures being the ones that find stronger evidence of change. The last line in the table, corresponds to comparing, for a fixed neighbourhood and year, working versus non-working days. Here, all methods but Wilcoxon's, find evidence of difference, with Schilling's method getting the strongest evidence.
Comparisons not included in this table, show that NO$_2$ contamination in working and non-working days is clearly different for all neighbourhoods and both years included in the data.

\subsection{Conclusion}

In the present article, we have discussed different ways of implementing the idea of permutation tests in the context
of the two-sample problem for functional data. The various approaches considered vary significantly in terms of
computational cost and in terms of power in different situations. Our discussion shows that
even the way in which a $p$-value
is assigned to a pair of statistics can be subject to significant variations.
At least two of the new methods proposed,
MA2 and the adaptation of Schilling's statistic, 
are highly competitive in terms of power against a broad range of alternatives,
as illustrated in the simulated and real data examples in Section 3. Thus, we expect to
have demonstrated sufficiently, the potential of the permutation test methodology 
in the context of functional data, and would expect this option to be 
considered and chosen in practical applications.

\begin{table}[H]
\caption{Empirical power (in \%)  of tests for geometric Brownian motion data}
\begin{center}
\begin{tabular}{|l|rrrrrr|}
\hline
Sample& HK& Wilcoxon &MA1&MA2&Schilling5& Schilling10 \\ \hline
Yx1.25 	&40	&4	&32	&45	&50	&67\\
Yx1.50 	&99	&3	&89	&100&100&100\\
Yx1.75 	&100&10	&100&100&100&100\\
Yx2.00 	&100&10	&100&100&100&100\\
Yr1.25	&49	&3	&15	&16	&31	&41\\    
Yr1.50 	&99	&4	&70	&86&97	&99\\
Yr1.75 	&100&8	&100&100&100&100\\
Yr2.00 	&100&11	&100&100&100&100\\
Ys1.25 	&9	&100&43	&98	&100&99\\
Ys1.50 	&21	&100&100&100&100&100\\
Ys1.75 	&33	&100&100&100&100&100\\
Ys2.00 	&39	&100&100&100&100&100\\
X 	   	&4	&4	&3	&3	&7	&1    \\ \hline
\end{tabular}
\end{center}
\label{tests}
\end{table}	

\begin{table}[H]
\caption{Hypothesis testing comparison on contamination data}\label{compcont}
\begin{center}
\begin{tabular}{lrrrr}
\hline
						& Wilcoxon & MA2 & Schilling ($k$ = 10) & HK \\ \hline
S2014n-w vs Pal2014n-w & 1.37E-04 & 0.002 & 0     &  0.06 \\
S2014n-w vs Pob2014n-w & 1.40E-05 & 0.001 & 0.137 &  0.05 \\
S2015w    vs E2015w      & 1.02E-04 & 0     & 0     &  0.001 \\
S2014n-w vs S2015n-w   & 0.385    & 0.576 & 0.62  &  0.07 \\
S2014w    vs S2015w      & 9.47E-01 & 0.002 & 0     &  0.001 \\ 
S2015w    vs S2015n-w   & 2.27E-03 & 0.008 & 0     &  0.001 \\ \hline
\end{tabular}
\end{center}
\label{testk10}
\end{table}

\end{document}